\documentclass[graybox]{svmult}

\usepackage{mathptmx}       
\usepackage{helvet}         
\usepackage{courier}        
\usepackage{type1cm}        
\usepackage{graphicx}        
\usepackage{multicol}        
\usepackage[bottom]{footmisc}

\begin{document}

\title*{Finite Elements for a  Beam System With Nonlinear Contact Under Periodic Excitation}
\author{H. Hazim, B. Rousselet}
\institute{H.Hazim \at Universit\'{e} de Nice Sophia-Antipolis, Laboratoire J.A. Dieudonn\'e U.M.R. C.N.R.S. 6621, Parc Valrose, F 06108 Nice, C{e}dex 2\, \email{hamad.hazim@unice.fr}
\and B. Rousselet\at Universit\'{e} de Nice Sophia-Antipolis, Laboratoire J.A. Dieudonn\'e U.M.R. C.N.R.S. 6621, Parc Valrose, F 06108 Nice, C{e}dex 2\, \email{br@math.unice.fr}}
%
%
\maketitle

\abstract{Solar arrays are structures which are connected to satellites; during
launch, they are in a folded position and submitted to high vibrations.
In order to save mass, the flexibility of the panels is not negligible
and they may strike each other; this may damage the structure. To
prevent this, rubber snubbers are mounted at well chosen points of the
structure; a prestress is applied to the snubber; but it is quite
difficult to check the amount of  prestress  and the snubber  may act
only on one side; they will be modeled as one sided springs (see figure $\ref{2}$).\\
In this article, some analysis for responses (displacements) in both time and frequency domains for a clamped-clamped Euler-Bernoulli beam model with a spring are presented. This spring  can be  unilateral or bilateral fixed at a point.
The  mounting (beam +spring)  is fixed on a rigid support which has a sinusoidal  motion of  constant
frequency.\\
The system is also studied in the frequency domain by sweeping frequencies between two fixed values, in order to save the maximum of displacements corresponding to each frequency.
Numerical results are compared with exact solutions in particular cases which already exist  in the literature.\\
On the other hand, a numerical and theoretical investigation of nonlinear normal mode (NNM) can be a new method to describe nonlinear behaviors, this work is in progress.\\
\begin{figure}[t]
\begin{center}
\includegraphics[scale=.25]{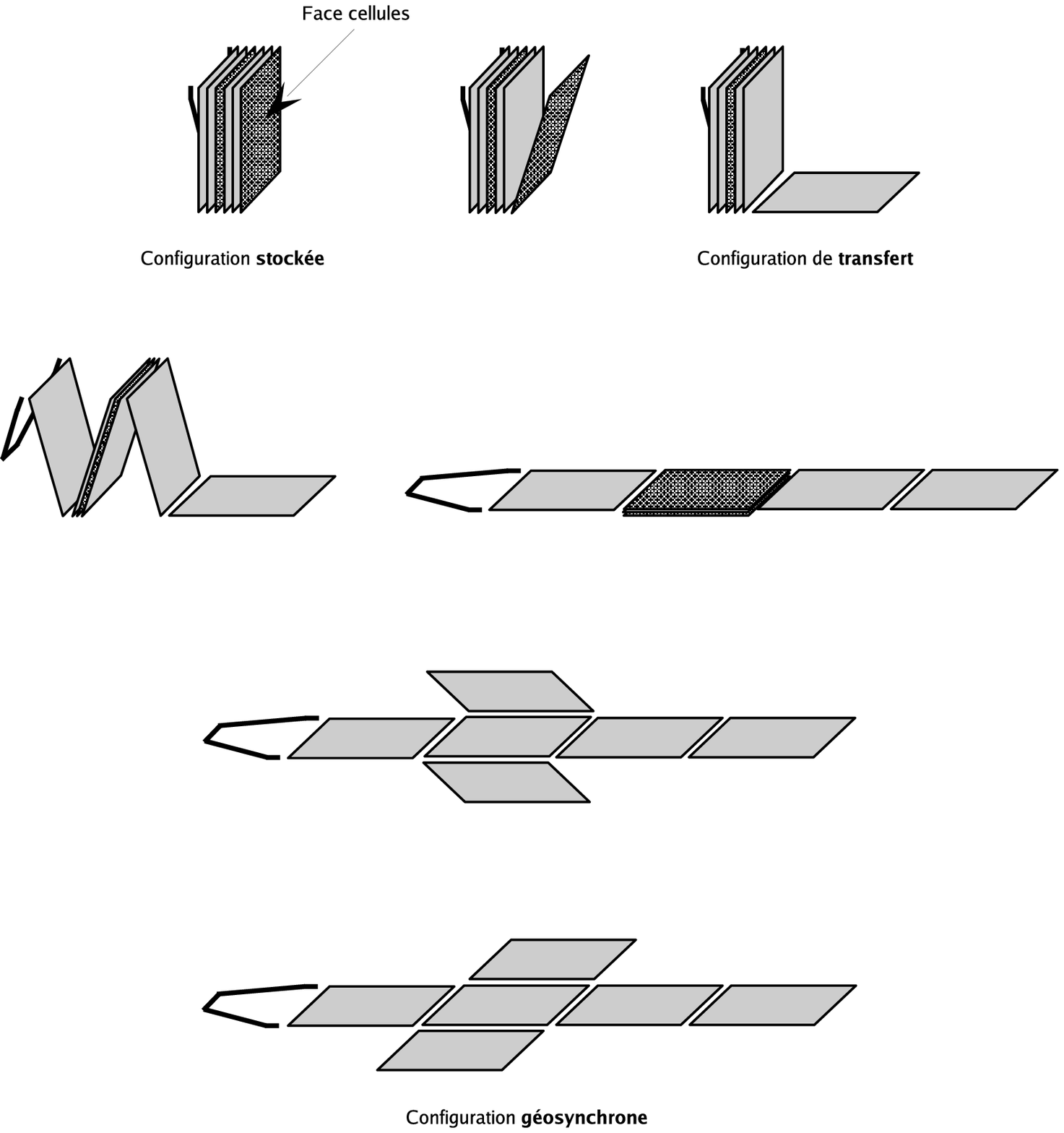}
\hspace{1.0cm}
\includegraphics[scale=.25]{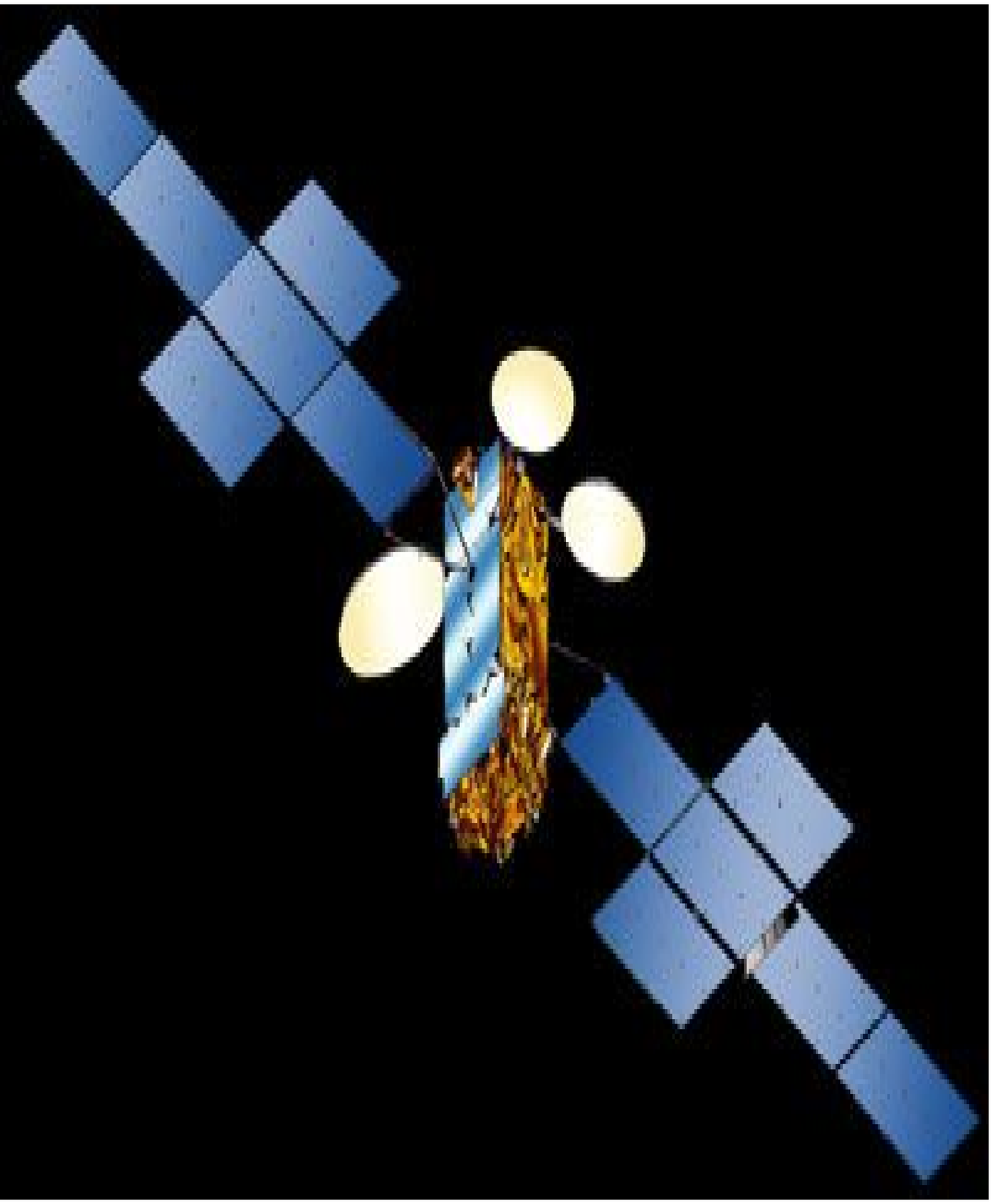}
\end{center}
\caption{\small{At left: solar arrays from folded to final position, at right: picture for the satellite AMC12 from Thales Alenia Space company}}
\label{1}
\end{figure}
}
\section{Introduction}
\label{sec:1}
\subsection{Previous Works}
In articles \cite{r1}, \cite{r2}, \cite{r4} et \cite{r5}, proposed by Thales Alenia Space research team, the dynamic of a beam sytem with a nonlinear contact force, under a periodic excitation given as an imposed acceleration form is studied both numerically and experimentally. When sweeping frequencies in an interval which contains eigen frequencies of the beam, resonance phenomena appear as well as new frequencies caused by the unilateral contact.\\
Finite element method in space domain is used, followed by numerical integration of the ordinary differential systems using specific software like 'STRDYN' of the finite element package DIANA.\\
The frequency sweeping  is done in different ways, one of these way is such as the  frequency $f$ changes as a function of the time $t$ according to :$f(t)=f_02^{st/60}$ where $s$ is the sweep rate in octaves/min and $f_0$ is the start frequency of the sweep. Results prove differences between sweep-down and sweep-up around eigen frequencies of the system where solutions are unstable.\\
At each value of the time $t$, computation is done and the maximum of displacement and acceleration are saved.\\Finally, comparison is made in time and frequency domains, between linear and nonlinear cases.
\subsection{Present Work}

This work is part of the phd work of the author under the guidance of B. Rousselet with the support of Thales Alenia Space, France.\\
Some analysis for responses (displacements) in both time and frequency domain for a clamped-clamped Euler-Bernoulli beam model with a linear spring are presented, the spring can be  unilateral or bilateral fixed at a point.
The  mounting (beam +spring)  is fixed on a rigid support which has a sinusoidal  motion of  constant
frequency.\\
The system is also studied  in the frequency domain by sweeping frequencies between two fixed values, such as saving the maximum of displacements corresponding to each frequency.

Numerical results are compared with exact solutions in particular cases which already exist  in the literature.\\
On the other hand, a numerical and theoretical investigation of nonlinear normal mode (NNM) can be a new method to distinguish  linear from nonlinear cases.\\

\section{Simplified Mechanical Model}
\label{sec:2}
\begin{figure}[t]
\includegraphics[scale=.25]{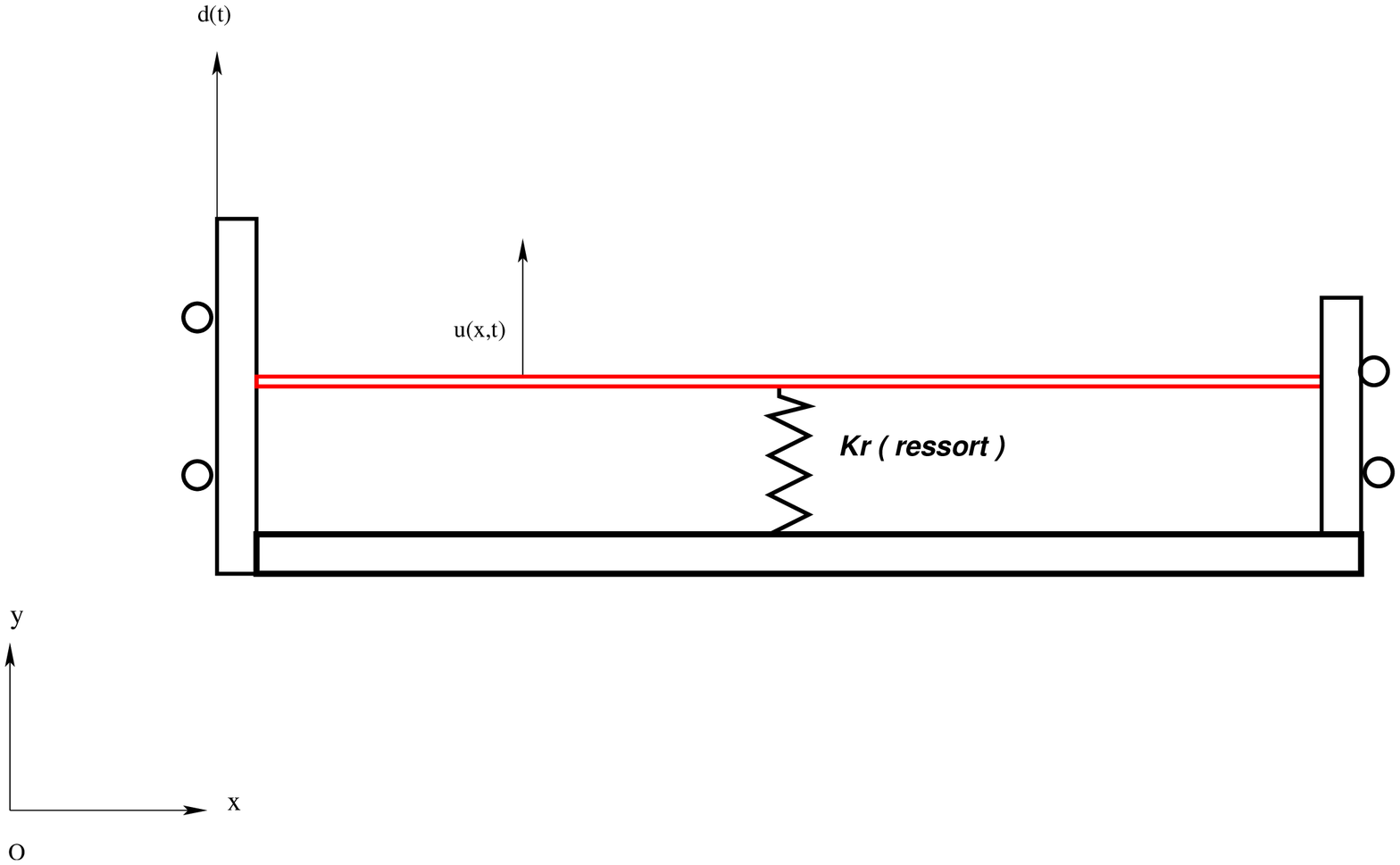}
\hspace{1cm}
\includegraphics[scale=.20]{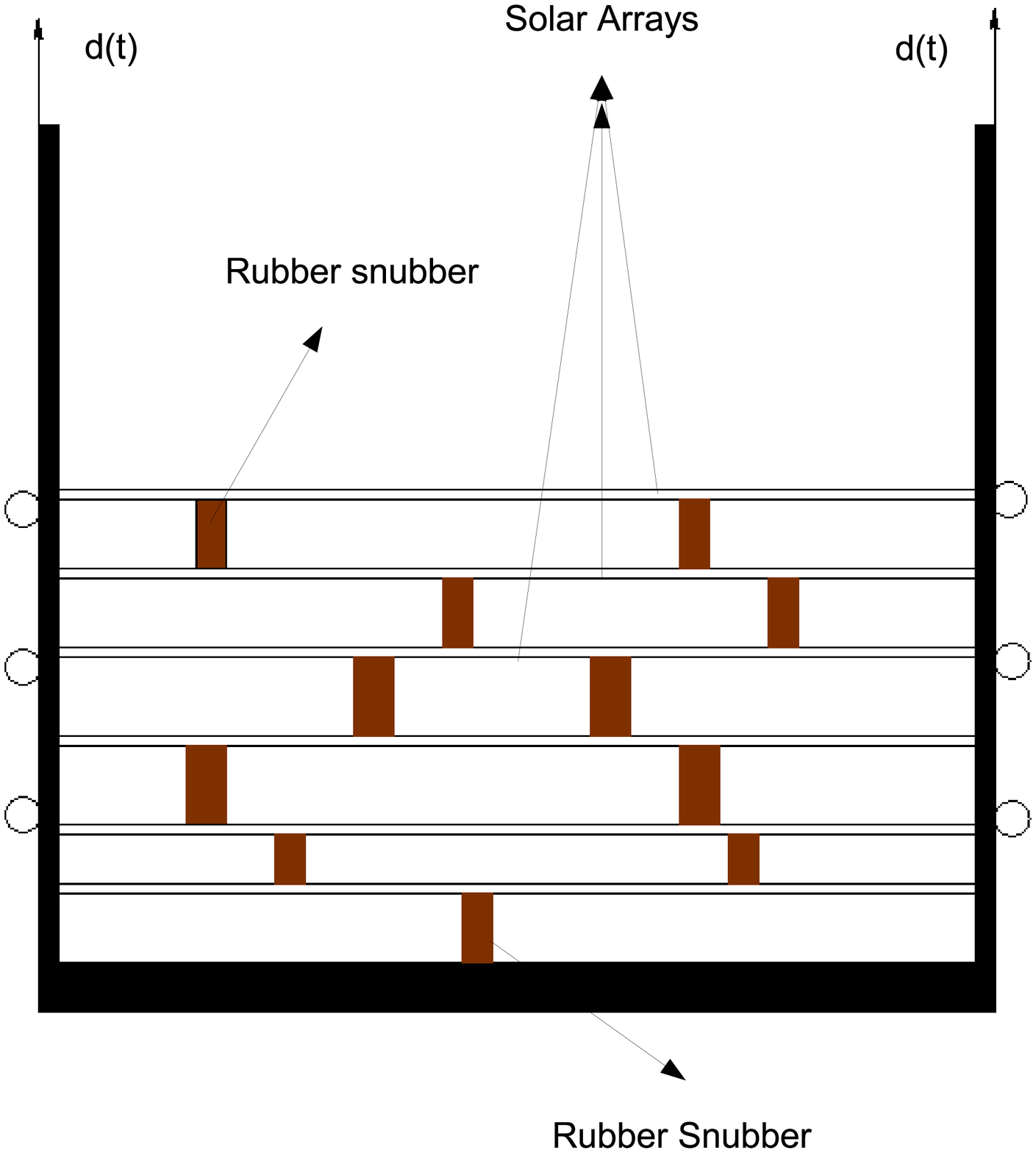}
\caption{ left: Simplified Mechanical Model, right: folded solar arrays with snubbers   }
\label{2}
\end{figure}

The study of the total  dynamic behavior of solar arrays in a folded position with snubbers  are so complicated, that to simplify, a solar array is modeled by a clamped-clamped Bernoulli beam with one-sided linear spring. This system is fixed on a shaker which has a vibratory motion $d(t)$ see figure ($\ref{2}$).

The motion of this beam system is modeled by the following PDEs with boundary conditions : $u(0,t)=u(L,t)=d(t)$ and $\partial_xu(0,t)=\partial_xu(L,t)=0$:\\ \\

\begin{equation}
\textbf{Bilateral}\;\textbf{ spring:}\;\;\rho S\ddot{u}(x,t)+EIu^{(iv)}(x,t)=k_{r}(d(t)-u(x_0,t))\delta_{x_0}\label{e1}
\end{equation}
\begin{equation}
\textbf{Unilateral}\; \textbf{spring:}\;\;\rho S\ddot{u}(x,t)+EIu^{(iv)}(x,t)=k_{r}(d(t)-u(x_0,t))_{+}\;\delta_{x_0} \label{e2}
\end{equation}
 \\
$u_+$ is the function defined by   $u_+=\frac{u+|u|}{2}$  \hspace{1cm} $L$= 0.485 m beam length. \\
$k_{r}$ = spring stiffness.  \hspace{2.0cm} $d(t)=-\frac{a}{(2\pi f)^2}\sin(2\pi ft)$ is the shaker motion.  \\
$\rho=2700$ $ kg/m^3$ beam density.  \hspace{2.0cm}  $S= 7,5.10^{-4} $ $m^2$ cross sectional area.\\
$E= 7.10^{10}$ $N/m^2$ Young's modulus. \hspace{0.5cm} $I= 1,41.10^{-8}$ $m^4$  second moment of area. \\
This data are  taken from $\cite{r2}$.
The classical Hermite cubic finite element approximation is used here to find an approximate  solution for equations $(\ref{e1})$ and (\ref{e2}), we find then two ordinary differential systems in the form :
\begin{equation}
M\ddot{q}+Kq=k_{r}(d(t)-q_{x_0})\overrightarrow{e_{x_0}} \label{s1}
\end{equation}
\begin{equation}
M\ddot{q}+Kq=k_{r}(d(t)-q_{x_0})_+\overrightarrow{e_{x_0}} \label{s2}
\end{equation}
$M$ et $K$ are respectively the mass and the stiffness assembled matrices, $q$ is the vector of  degree of freedom of the beam, $q_i=(u_i,\partial_x u_i)$, $i=1,2,...,n$. To each node are associated two degrees of freedom, the displacement and its  derivative.\\

To integrate numerically systems $(\ref{s1})$ and $(\ref{s2})$, we use the Scilab routines "ODE's" followed by the  FFT (Fast Fourier Transformation) to find frequencies of solutions, there is no special treatment for "ODE" routine to deal with the local non-differentiable  nonlinearity $(d(t)-u(x_0,t))_{+}$.
\section{Numerical Results}
\subsection{Highlights on the linearized system and the nonlinear effects }
\label{subsec:2}
There are two linear cases, the first one when there is no spring attached to the beam, and the second case is such as a linear spring is always attached to the beam ( bilateral spring ). The system becomes nonlinear when the spring becomes one-sided when the prestress is not well tuned. The nonlinearity has a special form, it is locally not differentiable but it is lipchitz.\\
Without spring, the motion has the following linear system of equation :
\begin{equation}
M\ddot{q}+Kq=0 \label{}
\end{equation}
The eigen frequencies of the motion are easily calculated  by computing the generalized eigen values of $M$ and $K$, this calculus can be done using software like  Scilab $\cite{a6}$. We are interested by the first three eigen frequencies of the system, we note that the precision depends on the number of the  finite elements used for the modeling of the system. We verify that ten finite elements give a good approximation for the first three  eigen frequencies  and their values are respectively :$334.21622Hz$, $921.48815Hz$ and $1807.7966Hz$. As known, the eigen-frequencies of a Bernoulli beam could be calculated using this formula:\\ $f_i=\frac{1}{2\pi}\sqrt{\mu_i^4\frac{EI}{Ml^4}}$,
$\mu_i$ $i=1,2,3,...n $ are given in $\cite{G}$. $\mu_1=4.73$, $\mu_2=7.853$ and  $\mu_3=10.996$ give $f_1=334.19889Hz$, $f_2=921.19996 Hz$ and $f_3=1806.1432 Hz$.
The nonlinearity on the system modifies the motion by adding new frequencies, subharmonics and superharmonics, besides  the eigen frequencies.
When the system is under periodic excitation, new frequencies appear also, there are many combinations of the excitation frequency with the system frequencies. These new frequencies appears in FFT of the system and also in the sweep test.
There are many ways to calculate  these frequencies, harmonic balance method and nonlinear normal mode (MNN) and asymptotic expansions methods   $(\cite{JB})$ .

\subsection{One Node Finite Element Model Without Periodic Excitation}
\label{subsec:2}
In this case, the beam is modeled by two finite elements without periodic excitation. Equations of displacement and its derivative are independent here because of the structure of   the mass and stiffness matrices:\\
 $M=\left(
                                                                                              \begin{array}{cc}
                                                                                                0.3647893 & 0 \\
                                                                                                0 & 0.0005500\\
                                                                                              \end{array}
                                                                                            \right)$ and $K=\left(
                                                                                                             \begin{array}{cc}
                                                                                                               1661090 & 0 \\
                                                                                                               0& 32560.825 \\
                                                                                                              \end{array}
                                                                                                            \right)$.\\
To find eigen  frequencies of the system in the linear case (system without contact), we juste have to calculate the generalized eigen values of $M$ and $K$.\\
The motion is divided into two phases: the first one when the beam touches  the spring and another one when the beam does not touch the spring.
The spring mass is negligible beside the beam mass, the equations of displacement of these two phases are respectively :\\
 $$M(1,1)\ddot{u} +(K(1,1)+k_r)u=0$$ $$M(1,1)\ddot{u}+K(1,1) u =0$$

\begin{figure}[t]
\begin{center}
\includegraphics[scale=.65]{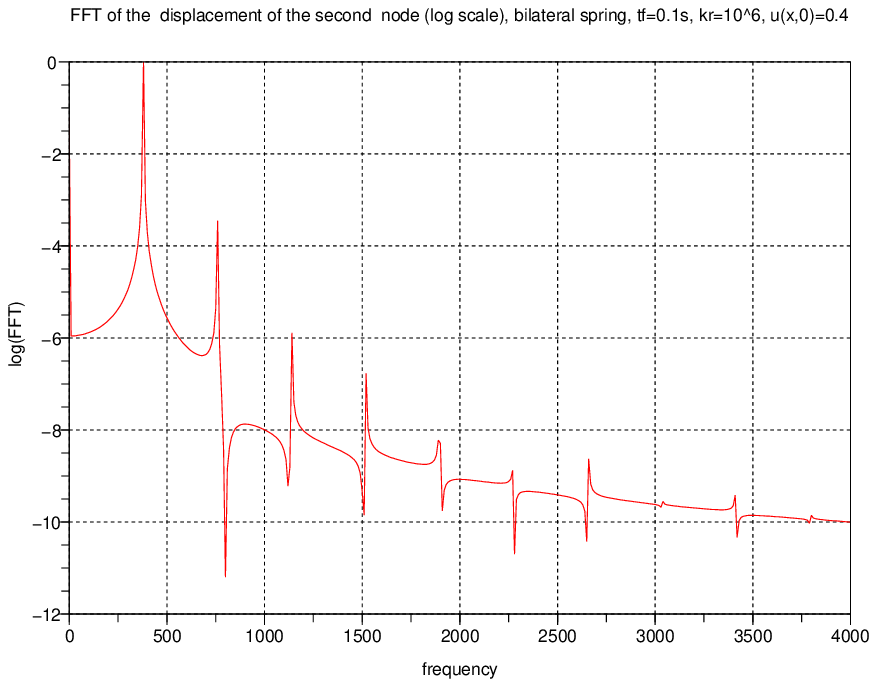}
\includegraphics[scale=.20]{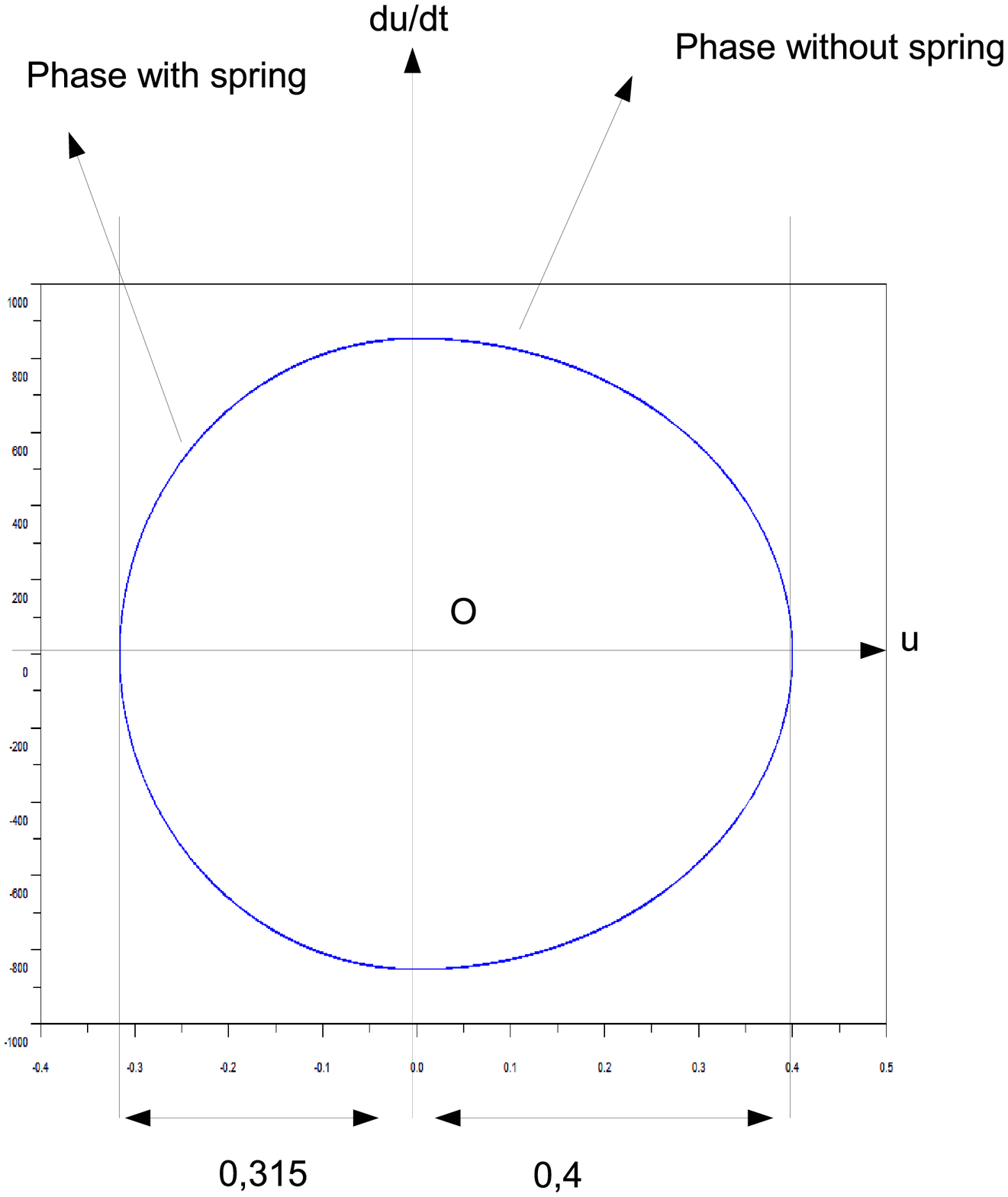}
\end{center}
\caption{\small{ Left: FFT  of the displacement of the second node with unilateral spring contact with initial conditions $u(x,0)=0.4$, $\dot{u}(x,0)=0$. Right: solution in the phase space.  }}
\label{3}
\end{figure}

The boundary conditions are :$u(0,t)=u(L,t)=0 $ and $\partial_xu(0,t)=\partial_xu(L,t)=0$.\\
The period  of the solution will be the sum of the half period of the first phase and the second phase :
$$T=\frac{\pi}{\sqrt{\frac{K_{11}(1,1)}{M_{11}(1,1)}}}+\frac{\pi}{\sqrt{\frac {K_{11}(1,1)+k_r}{M_{11}(1,1)}}}$$
the numerical value of the motion frequency $\frac{1}{T}$ is $384.74186 Hz$ this is in agreement with the numerical calculus.
In figure ($\ref{3})$, The Fast Fourier Transformation shows frequencies of the system, the first peak correspond to the analytical value. The motion in the phase plane shows a periodic conservative solution.


\subsection{One Node  Finite Element Model, Beam Under Periodic Excitation } \label{1}
\label{subsec:2}
The beam is modeled by two finite elements under periodic prescribed displacement $d(t)$. The spring is always fixed in the middle of the beam on a node, it can be unilateral or bilateral.\\

\begin{figure}[]
\includegraphics[scale=.65]{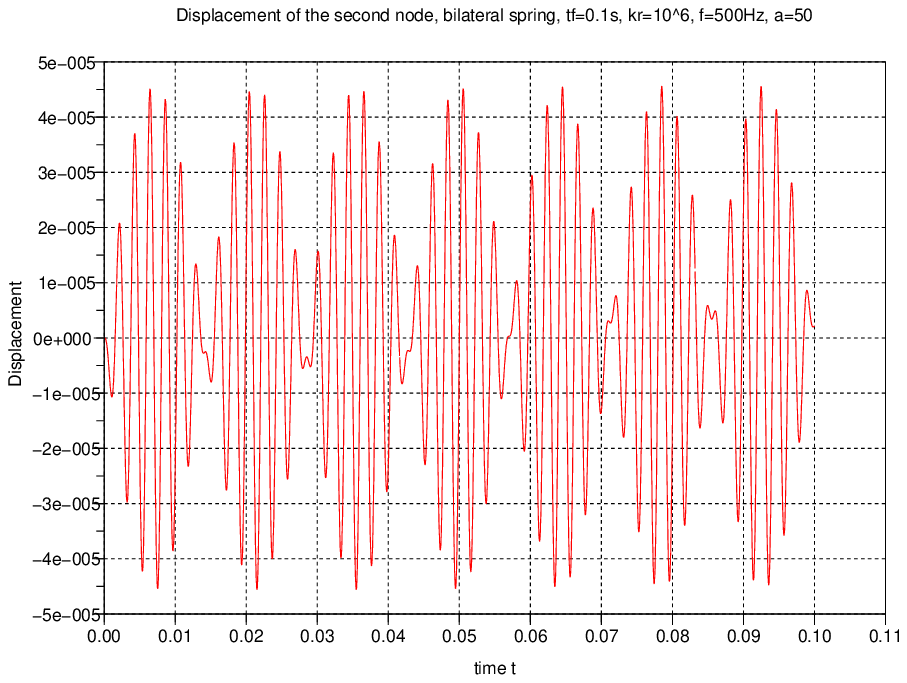}
\includegraphics[scale=.65]{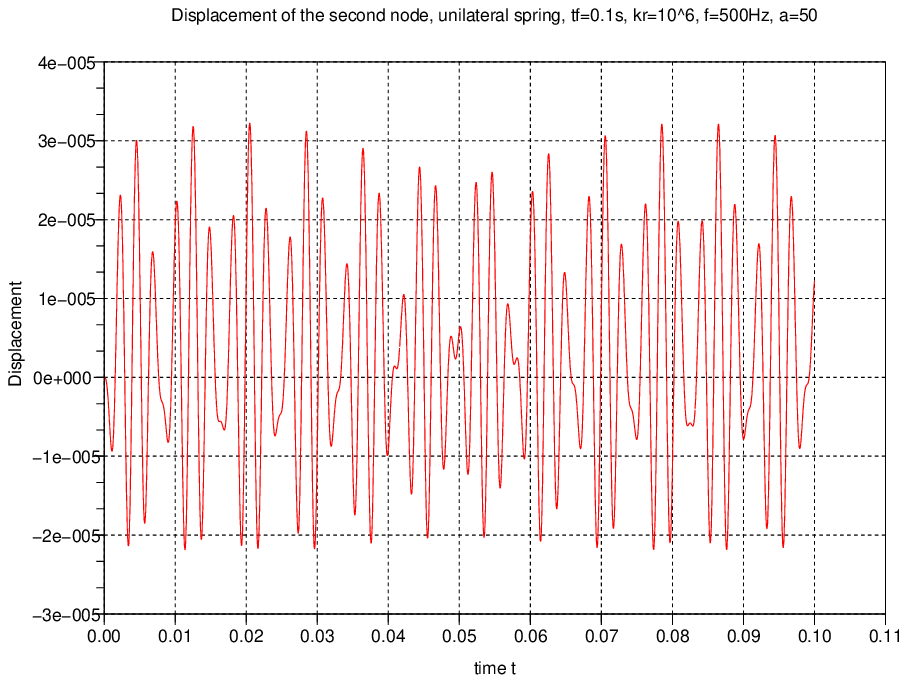}
\caption{ Displacement of the second node with \textbf{bilateral and unilateral }spring contact under periodic excitation of 500 Hz (Two finite element model) }
\label{4}
\end{figure}

Figure ($\ref{4}$) shows the displacement of the second node (in the middle) of the beam with bilateral and unilateral spring, the period of the motion changes with unilateral spring. Figure ($\ref{5}$) shows frequencies of the beam with bilateral and unilateral cases; In the bilateral case, there are two frequencies, the eigen frequency (430 Hz)  and the frequency of the shaker (500 Hz); In the unilateral case, there are the system frequency, the shaker frequency and many other superharmonics and subharmonics frequencies due  to the contact.


\subsection{Ten Finite Elements}
\label{subsec:2}
In this section, the beam is modeled with ten finite elements, this approximation is quite good if we are interested by the first two eigen frequencies of the system.\\
The spring is always fixed in the middle of the beam on a node, and can be unilateral or bilateral, the whole system is fixed  on a shaker which has periodic prescribed displacement $d(t)$, we compare the linear with the nonlinear case for the same value of parameters.\\
Figure ($\ref{6}$) shows the displacement in the time domain for bilateral and unilateral case and their corresponding frequencies using FFT.
In the linear case (bilateral spring), the first and the third eigen frequencies are shown, the shaker frequency also(500 Hz); The middle of the beam is node of the second mode. In the nonlinear case (unilateral spring), frequencies corresponding to the first and the third eigen frequencies of the linear case and also the shaker frequency are shown; Many other superharmonics, subharmonics and combination with the shaker frequency appear also due to the non linear contact.

\section{Frequency Sweep Excitation}
Frequency sweep excitation is usually used in the experiments to check the dynamic behavior of mechanical systems.
As mentioned before, the beam is under periodic excitation given as a  displacement form\\ $d(t)=-\frac{a}{(2\pi f)^2}\sin(2\pi ft)$, that means that the beam is under an effort (acceleration)  $\ddot{d}(t)=a\sin(2\pi ft)$, $a$ in $m/s^2$ is the amplitude of the acceleration and $f$ in $Hz$ is the frequency of excitation.

We are just interested by a sweep-up test, for a initial given  frequency $f_0$ and a fixed value of $a$. We compute the solution of differential systems, then we save the maximum of the acceleration and the displacement. In a second time, we add a fixed value to $f_0$, the frequency step $df$.
We compute again for $f_0+df$ and we save the same quantities for this iteration, the initial conditions are fixed on zero again.
We continue our test to reach a fixed frequency $f_1$ chosen such as to cover the first and the second eigen frequency of the system.

Finally, we plot the maximum saved in each sweep test, then we compare linear with nonlinear cases, we study also the effect of the amplitude $a$ and the frequency step on the system.\\

We remark that obtained curves are similar to curves of $FFT$, it is an another method to find frequencies for  mechanical systems.

\subsection{Two Finite Elements}\label{two}
\label{subsec:2}
The beam is modeled by two finite elements, the spring is fixed in the middle, the Bernoulli beam is clamped in its  both extremities, we just have a free node of two degrees of freedom, the first is the displacement, the second is the derivative of the displacement.
The eigen frequency of the linearized system in this case is around  $339Hz$, the sweep-up  begin from $100Hz$ to $1000Hz$, the frequency step $df=5Hz$ and the amplitude of excitation $a=50m/s^2$.\\

\begin{figure}[]
\includegraphics[scale=.65]{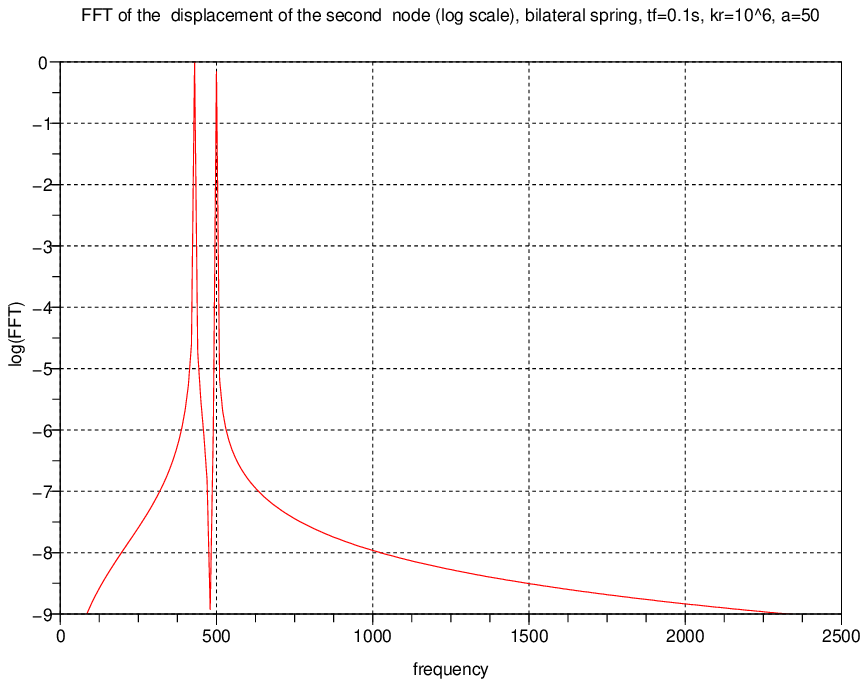}
\includegraphics[scale=.65]{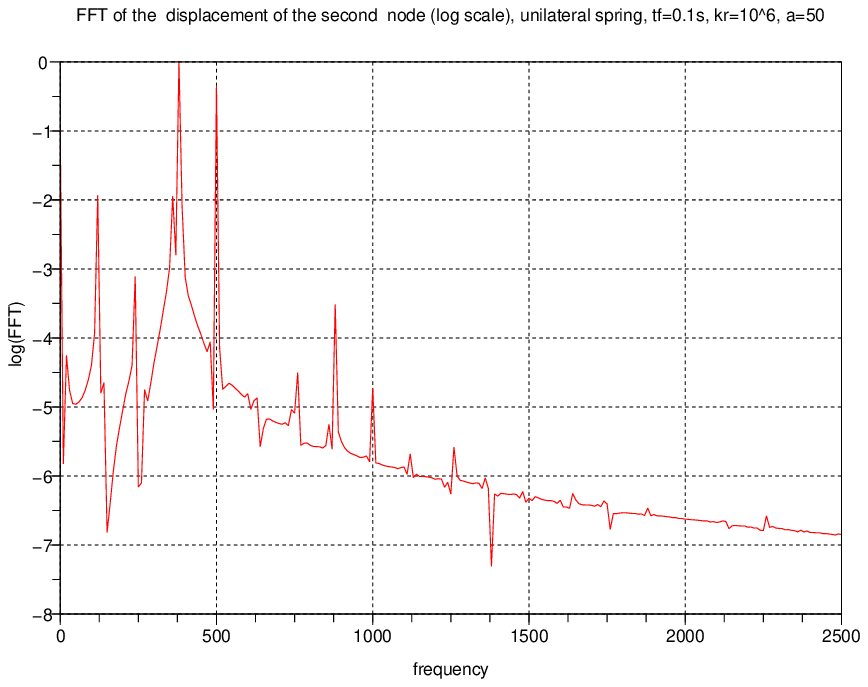}
\caption{ FFT of the displacement of the second node with \textbf{bilateral and unilateral }spring contact under periodic excitation of 500 Hz (Two finite element model)}
\label{5}
\end{figure}

\begin{figure}[]
\includegraphics[scale=.65]{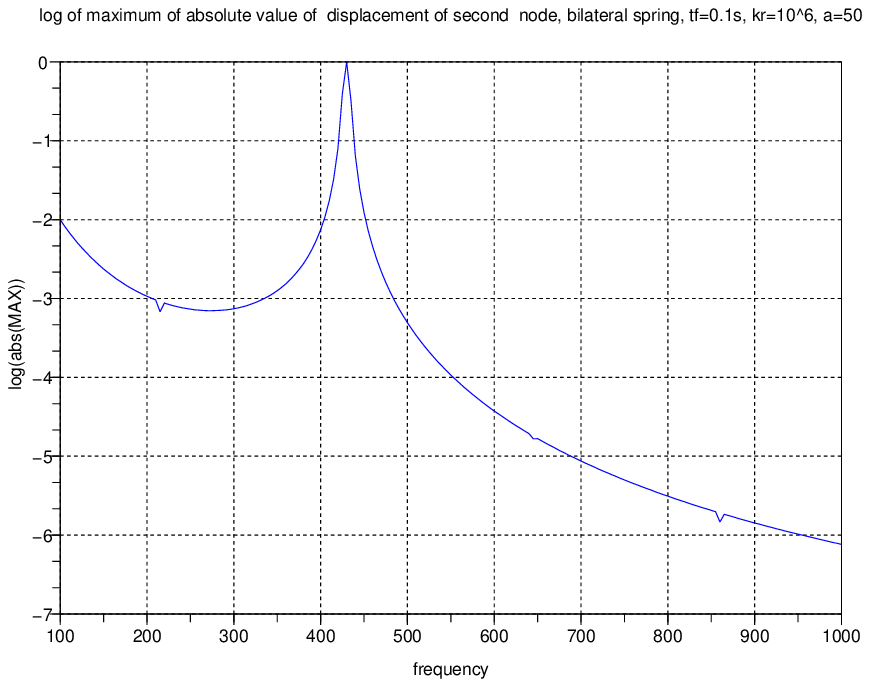}
\includegraphics[scale=.65]{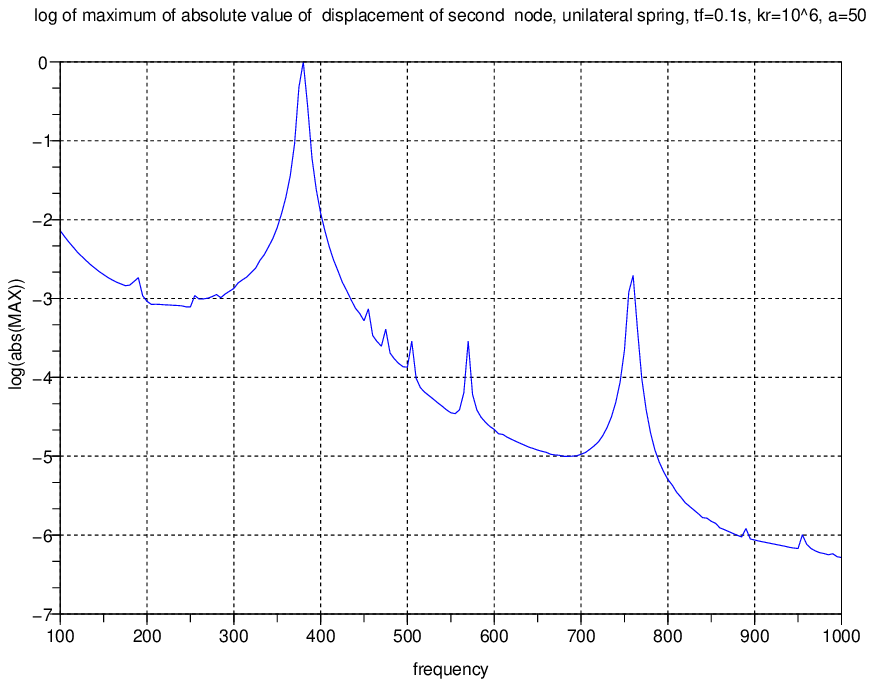}
\caption{ \small{Maximum of absolute value of the displacement  for both \textbf{bilateral and unilateral spring}, two finite element model, $kr=10^6 N/m$, $a=50m/s^2$, $tf=0.1s$} (Sweep-up test)}
\label{7}
\end{figure}

In each iteration, the integration time is $tf=0.1s$, the initial conditions are always fixed at $0$ : $q(x,0)=\dot{q}(x,0)=0$, $k_r=10^6N/m$.\\
In figure ($\ref{7}$), the peak in the bilateral case corresponds to the eigen frequency of the system, its  abscissa  is $430Hz$. In the unilateral case, the peak corresponds to the eigen frequency too, its abscissa is $384Hz$. Other peaks appear; they are due to the unilateral contact. These results conform well with the FFT in figure ($\ref{5}$).

\subsection{Ten finite elements}
\label{subsec:2}
\begin{figure}[]
\includegraphics[scale=.65]{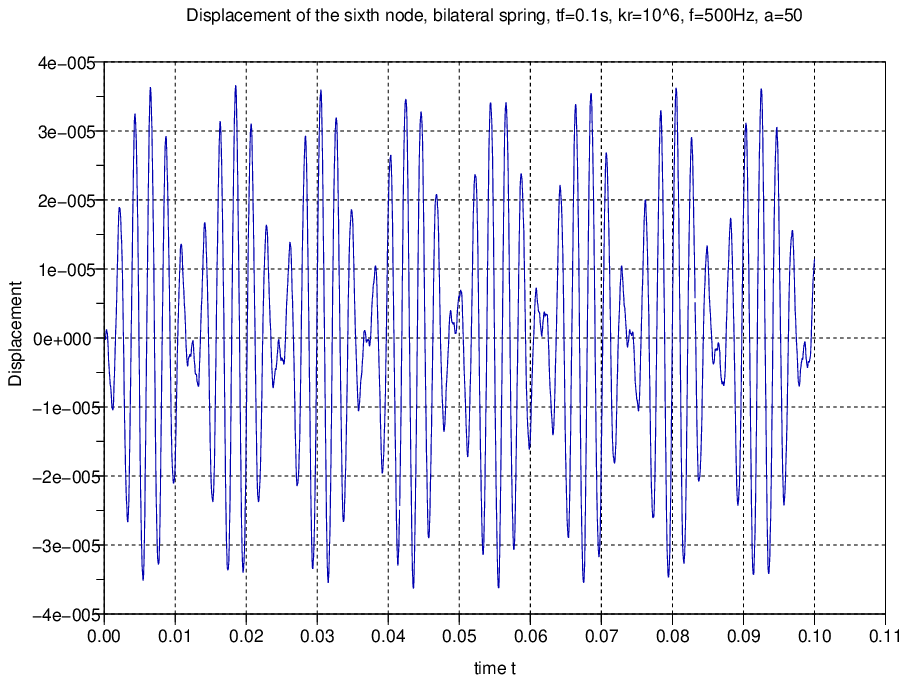}
\includegraphics[scale=.65]{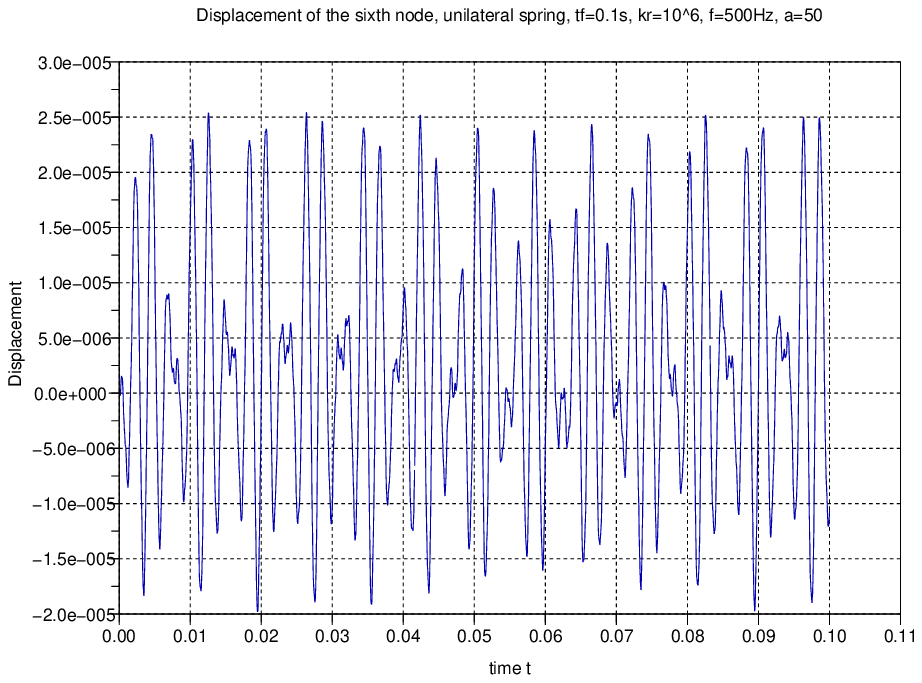}
\includegraphics[scale=.65]{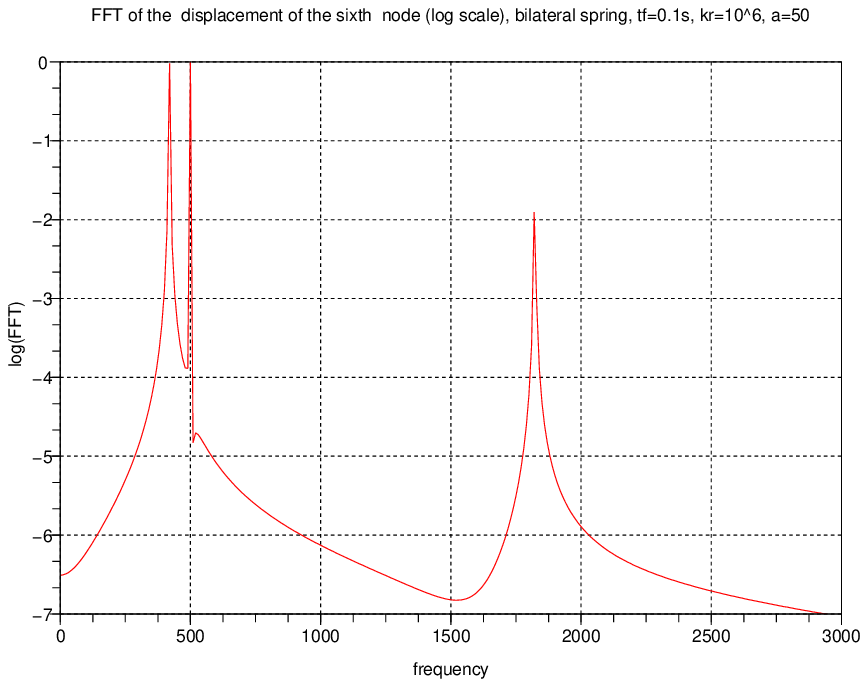}
\includegraphics[scale=.65]{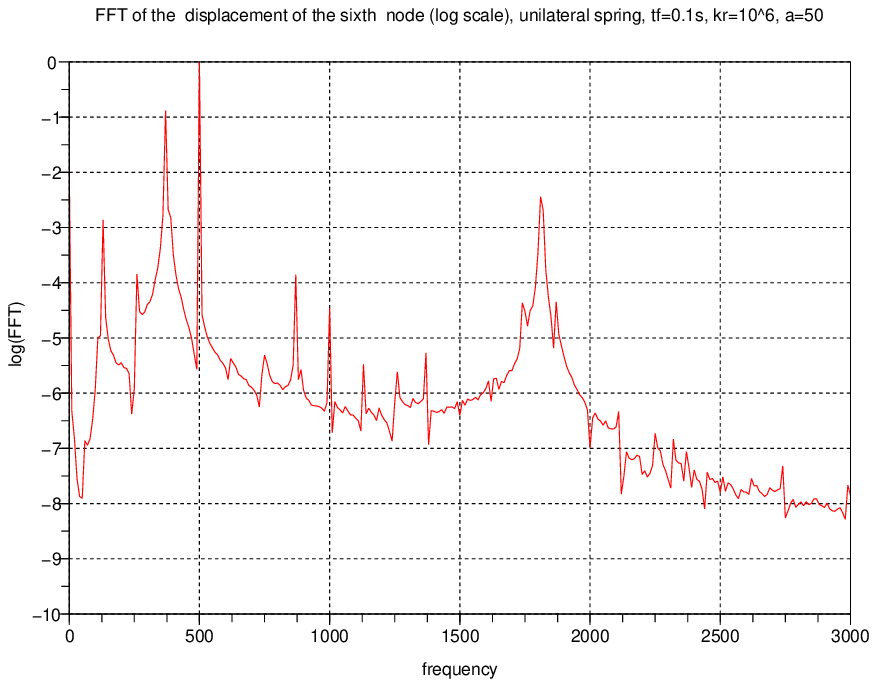}
\caption{ Displacement and their FFT of the sixth node with \textbf{bilateral and unilateral} spring contact under periodic excitation of 500 Hz (Ten finite element model) }
\label{6}
\end{figure}

\begin{figure}[]
\includegraphics[scale=.65]{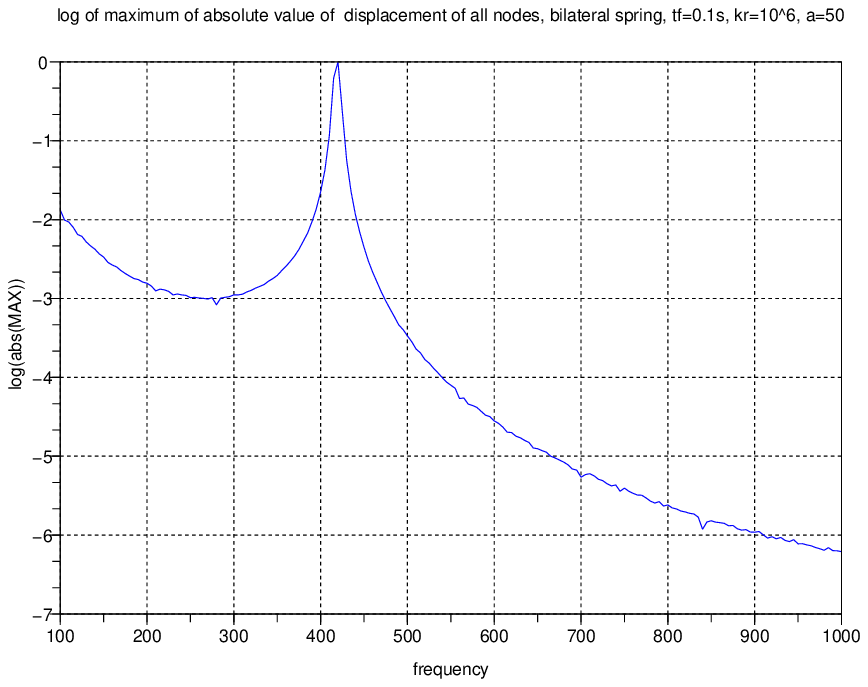}
\includegraphics[scale=.65]{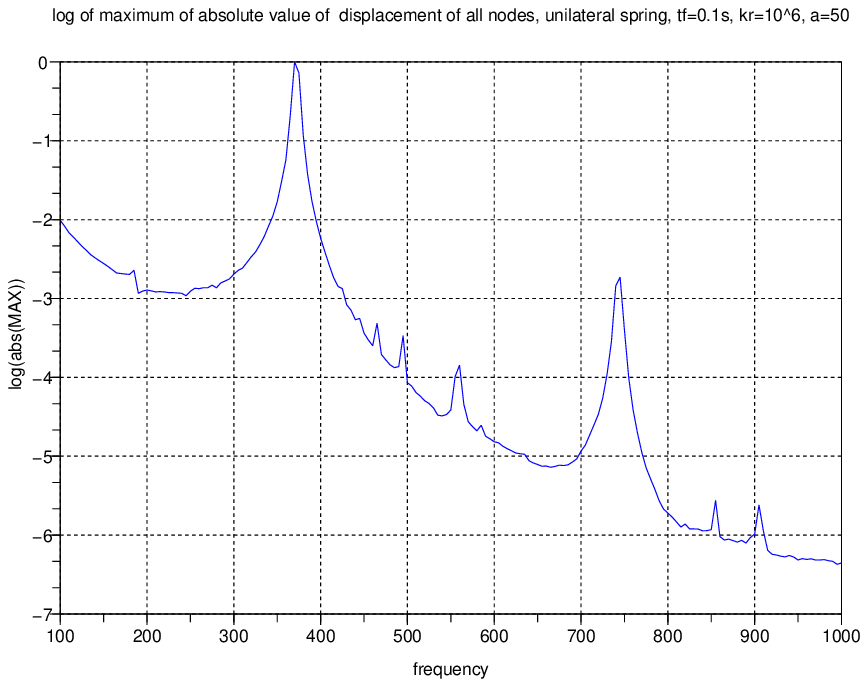}
\caption{\small{Maximum of absolute value of the displacement of all nodes(log scale) for both \textbf{bilateral and unilateral spring}, ten finite element model $k_r=10^6$, $tf=0.1s$, $a=50m/s^2$} (sweep-up test) }
\label{8}
\end{figure}

The beam is modeled by ten finite elements, the spring always in the middle, the other parameters are the same defined in section $\ref{two}$.\\
In figure ($\ref{8}$), the peak in the bilateral case  corresponds to the first eigen frequency of the system, the second eigen frequency does not appear because the imposed displacement is not enough to excite it ($d(t)=-\frac{a}{(2\pi f)^2}\sin(2\pi ft)$ very small when f become so high).
Peaks in the unilateral case show the system frequencies corresponding  to the eigen frequency of the linear system and other frequencies due to the contact. These results conform with the FFT in figure ($\ref{6}$) in the interval $[100,1000]$ Hz.

\newpage
\section{Conclusion}
We have presented some preliminary numerical results to compare the
vibrations of a beam equipped with a bilateral or a unilateral spring.
Asymptotic expansions using some results of S. Junca ($\cite{JB}$) are in
project and compared with numerical results in order to asses the
quality of both approaches. In particular  accurate computations of non
linear normal modes (see $\cite{a4}$)  will be considered.

\end{document}